\magnification=\magstep1 
\overfullrule=0pt      \input epsf    \voffset=.4in 
\def\eqde{\,{\buildrel \rm def \over =}\,} 
\def\leaderfill{\leaders\hbox to 1em{\hss.\hss}\hfill} 
 \def\la{\lambda}   
       \def\i{{\rm i}}

  \def\L{{\Lambda}} \def\M{{\cal M}}

\font\huge=cmr10 scaled \magstep2

\font\smcap=cmcsc10       
  \def\G{\Gamma}
   
\input amssym.def
\def\Z{{\Bbb Z}} \def\R{{\Bbb R}} \def\Q{{\Bbb Q}}  
\def\C{{\Bbb C}}   \def\M{{\Bbb M}}
\def\H{{\Bbb H}}

\centerline{{\bf \huge Postcards from the Edge,}}\medskip
\centerline{{\bf \huge or Snapshots of the Theory of Generalised
Moonshine}\footnote{$^{\dag}$}{{\sevenrm This is the text of my talk at
the Banff conference in honour of R.V. Moody's 60th birthday. A stream-lined version of this
paper (with the pedagogy removed) is my contribution to
a volume in his honour.}}} \bigskip
\bigskip   \centerline{Terry Gannon}
\medskip\centerline{{\it Department of Mathematical Sciences, University of Alberta,}}
\centerline{{\it Edmonton, Alberta, Canada, T6G 2G1} } \smallskip
\centerline{{e-mail: tgannon@math.ualberta.ca}}
\bigskip\bigskip

\rightline {\sevenrm I dedicate this paper to a man who throughout his
career}

\rightline{\sevenrm has exemplified the power of conceptual thought in math: Bob Moody.}\bigskip

In 1978, John McKay made an intriguing observation: $196\,884\approx
196\,883$. {\it Monstrous Moonshine} is the collection of questions
(and a few answers) inspired by this observation. In this paper we
provide a few snapshots of what we call the underlying theory. But
first we digress with a quick and elementary review.\medskip 

By a {\it lattice} in $\C$ we mean a discrete subgroup of $\C$ under
addition. We can always express this (nonuniquely) as the set of points
$\Z w+\Z z\eqde\L\{w,z\}$. We dismiss as too degenerate the lattice $\L=\{0\}$.
Call two lattices $\L,\L'$ {\it similar} if
they fall into each other once the plane $\C$ is rescaled and rotated 
about the origin --- i.e.\
$\L'=\alpha \L$ for some nonzero $\alpha\in\C$.  In Figure 1 we draw
(parts of) two similar lattices.  For another example, consider the degenerate
case where $w$ and $z$ are linearly dependent over $\R$: then in fact
$w$ and $z$ are linearly dependent over $\Z$
(otherwise discreteness would be lost) and any such lattice is
similar to $\Z\subset \C$.

\medskip\epsfysize=1.8in\centerline{ \epsffile{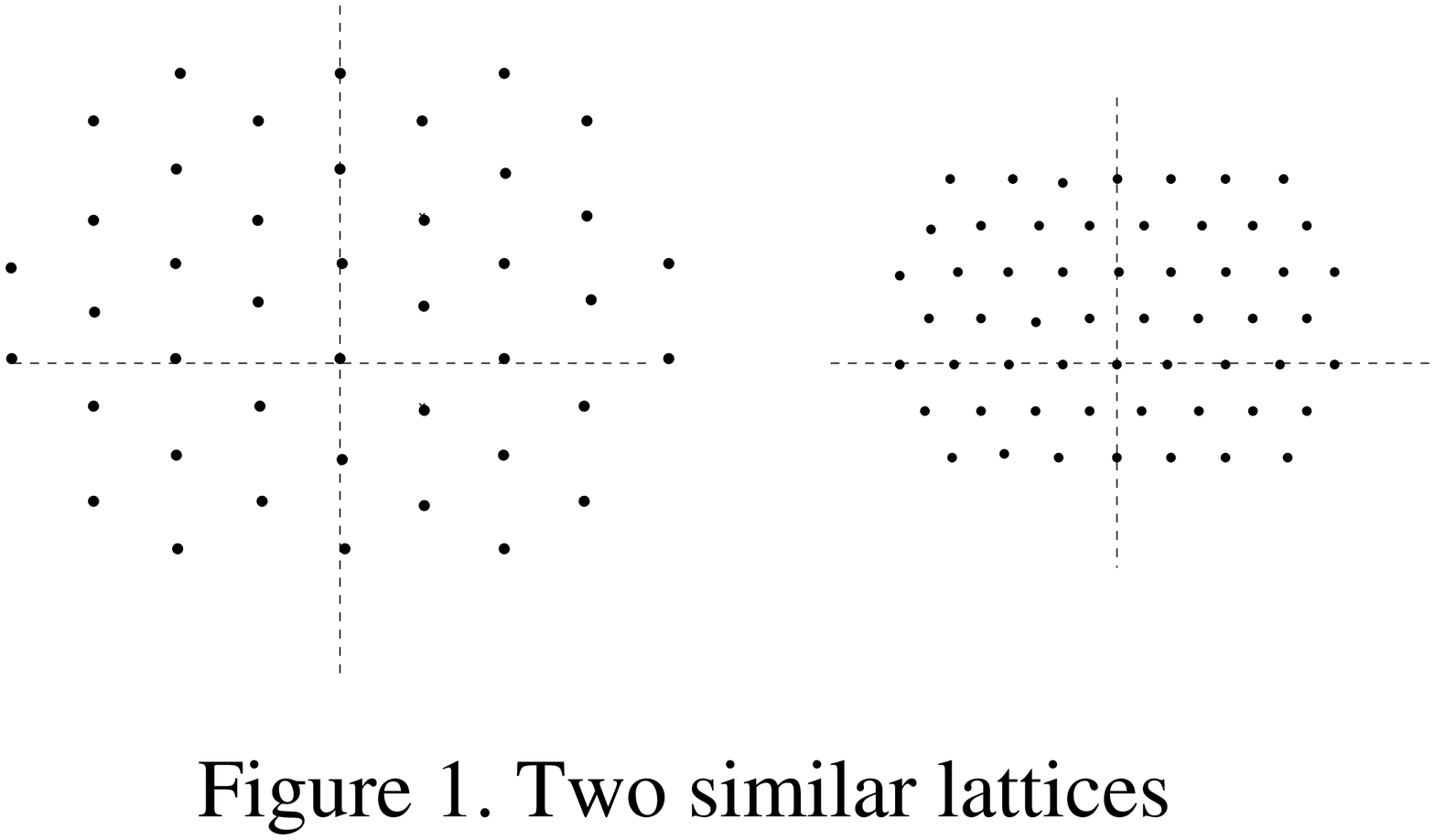}}\medskip

We're interested in the set of all equivalence classes $[\L]$ of
similar lattices. There is a natural topology on this set, and in fact
a differentiable structure. Now, it's a lesson of modern geometry that one
probes a topological set by considering the functions which live on
it. So consider  any complex-valued function $f(\L)$ on the
set of all lattices $\L$, which maps similar lattices $\L,\L'$ to the identical
number $f(\L)=f(\L')$ --- in other words, $f$ is well-defined on the equivalence classes
$[\L]$. It also is required to be a `meromorphic
function' of these classes $[\L]$. To specify what this means, we need
to look more closely at the set of all $[\L]$.

Let $w$ and $z$ be linearly independent over $\R$ (the
generic case). By choosing $\alpha=\pm 1/w$, we get that any lattice
$\L\{w,z\}$ is similar to one of the form $\L(\tau)\eqde\Z\tau+\Z$
where Im$(\tau)>0$. So each equivalence class $[\L]$ of (nondegenerate) similar lattices  can
be associated to a point $\tau$ in the {\it upper half-plane} $\H$.
In the same way, the degenerate ones are assigned
$\tau\in\Q\cup\{\infty\}$. In other words we can regard our 
function $f([\L])$ as being a complex-valued function on
$\overline{\H}\eqde\H\cup\Q\cup\{\infty\}$. We require it to be
meromorphic (i.e.\ complex-analytic apart from isolated poles) on $\H$, and
also meromorphic at $ \Q\cup\{\infty\}$ (we'll define this shortly).
Since
$\L\{w,z\}=\L\{z,w\}$, we get that the lattices
$\L(\tau)$ and $\L(-1/\tau)$ are similar. In fact more
generally the equivalence classes $[\L(A.\tau)]=[\L(\tau)]$ are equal,
for any matrix $A\in{\rm SL}_2(\Z)$ --- the group SL$_2(\Z)$ consists of all
$2\times 2$ integer matrices with determinant $\pm 1$, and
$A=\left(\matrix{a&b\cr c&d}\right)$ acts on
$\tau\in\overline{\H}$ by the fractional
linear transformation $A.\tau={a\tau+b\over c\tau+d}$.
Our function $f:\overline{\H}\rightarrow\C$ thus
must have the group SL$_2(\Z)$ as its symmetry group:
$f(A.\tau)=f(\tau)$ for all matrices $A\in{\rm SL}_2(\Z)$. In fact, this
is the only redundancy in our identification of equivalence classes of
similar lattices with points in $\overline{\H}$: each class $[\L]$ corresponds to 
precisely one SL$_2(\Z)$-orbit in $\overline{\H}$. This space
SL$_2(\Z)\backslash\overline{\H}$ of orbits is called the {\it moduli space}
for similar lattices, the simply-connected space $\overline{\H}$ is
called its {\it Teichm\"uller} or {\it universal covering space}, and the
redundancy group SL$_2(\Z)$ (or really PSL$_2(\Z)$) is called its
{\it modular group}.

We still have to explain what we mean by `meromorphic at the
degenerate lattice class $[\Z]$', i.e.\ at $\tau\in\Q\cup\{\infty\}$. It's enough to
consider $\tau=\infty$, by the SL$_2(\Z)$ symmetry. Because $f(\tau)$
has period 1, we can expand $f(\tau)$ as a power series in the variable
$q=e^{2\pi\i\tau}$: $f(\tau)=\sum_{n\in\Z}a_nq^n$. We require that
$a_n=0$ for all $n$ sufficiently close to $-\infty$ --- in other words the only
possible singularity at $q=0$ is a pole of finite order. These degenerate
points $\tau\in\Q\cup\{\infty\}$ are called cusps.

\medskip\noindent{{\smcap Definition 1.}} A {\it modular function} $f$
(for SL$_2(\Z)$) is a meromorphic function $f:\overline{\H}\rightarrow\C$,
obeying the symmetry
$f(A.\tau)=f(\tau)$ for all $\tau\in\overline{\H}$ and $A\in{\rm SL}_2(\Z)$.\medskip

We can construct some modular functions as follows. Define the {\it
(classical) Eisenstein series} by
$$G_k(\tau)\eqde\sum_{{m,n\in\Z\atop (m,n)\ne (0,0)}}(m\tau+n)^{-k}\eqno(1)$$
For odd $k$ it identically vanishes. For even $k>2$ it converges
absolutely, and so defines a function holomorphic throughout $\H$. It
also is holomorphic at the cusp $\tau=\infty$. We get 
$$G_k({a\tau+b\over
c\tau+d})=(c\tau+d)^k\,G_k(\tau)\qquad\forall\left(\matrix{a&b\cr
c&d}\right)\in{\rm SL}_2(\Z)\eqno(2)$$
and all $\tau$. This is because the sum in (1) is really over all
nonzero $x\in\L(\tau)$, and because SL$_2(\Z)$  parametrises certain
changes-of-basis $\{\tau,1\}\mapsto\{w,z\}$ of the two-dimensional lattice $\L(\tau)$. This
transformation law (2) means that $G_k$ isn't quite a modular function
(it's called a {\it modular form}). However,
various homogeneous rational functions of these $G_k$ {\it will} be
modular functions --- for example $G_8(\tau)/G_4(\tau)^2$ (which turns
out to be constant) and $G_4(\tau)^3/G_6(\tau)^2$ (which doesn't).
We'll see shortly that all modular functions arise in this way.

Why are modular functions interesting? At least in part, this has to
do with the omnipresence of two-dimensional lattices. For instance,
given a (nondegenerate) lattice $\L=\L\{w,z\}$, the quotient $\C/\L$
is a torus; the converse is true too. Similar lattices correspond to conformally
equivalent tori (i.e.\ the equivalence preserves the angles between
intersecting curves, but not
necessarily arc-lengths). So a modular function lives on the moduli space
of conformally equivalent tori. To push this thought a little further,
the circle $x^2+y^2=1$ is
really a sphere (with two points removed) embedded in $\C^2$ when $x,y$
are regarded as
complex variables (to see this, consider the change of variables
$x={1\over 2}(w+w^{-1}),y={\i\over 2}(w-w^{-1})$). In the same way,
the (nondegenerate) cubic $y^2=x^3+ax^2+bx+c$ is really a torus (with
one point removed) in $\C^2$ when $x,y$ are complex --- the torus
$\C/\L(\tau)$ corresponds to the complex curve $y^2=4x^3-ax-b$ where
$a=60\,G_4(\tau)$ and $b=140\,G_6(\tau)$.  (Incidentally, the missing
points on the sphere and torus appear naturally when projective
coordinates are used, i.e.\
the missing points are `points at infinity'). More precisely, (birational)
equivalence classes of cubics (over $\C$) are parametrised by our familiar
moduli space SL$_2(\Z)\backslash\H$ --- conformal structure on a real
surface corresponds to complex-differentiable structure on the corresponding
complex curve. So modular functions
can arise whenever tori or cubics (more commonly called `elliptic curves') arise.
Elliptic curves are special because they're the only complex
projective curves which have an algebraic group structure.
In any case, modular functions and their various generalisations hold a 
central position in both classical and modern number theory. For a very enjoyable account
of  the classical theory, see [21].

Can we characterise all modular functions? The key idea is to look directly at
the moduli space $M={\rm SL}_2(\Z)\backslash\overline{\H}$. We know that any modular
function will be a meromorphic function on the surface $M$. Thanks to the presence
of the cusps (i.e.\ the degenerate lattices), $M$ will be a compact Riemann
surface. With a little bit of work, it can be quite easily seen that
it is in fact a sphere. Although there are large numbers of
meromorphic functions on the complex plane $\C$, the only ones of
these which are also meromorphic at $\infty$ are the rational
functions ${{\rm polynomial\ in}\ z\over {\rm polynomial\ in}\ z}$ (the
others have essential singularities there). In
other words, the only functions meromorphic on the Riemann sphere
$\C\cup\{\infty\}$ are the rational functions. So if $J$ is a
change-of-variables (or {\it uniformising}) function from our moduli space $M$ to the Riemann
sphere, then $J$ (interpreted as a function on the covering space
$\overline{\H}$) will be a modular function, and any modular function $f(\tau)$ will be a
rational function in $J(\tau)$: $f(\tau)= {{\rm polynomial\ in}\ J(\tau)\over {\rm polynomial\ in}\ J(\tau)}$.
And conversely, any rational function in $J$ will be modular.
Thus $J$ generates modular functions, in a way analogous to (but stronger
and simpler than) how the exponential $e(x)=e^{2\pi\i\,x}$ generates the
period-1 continuous functions $f$ on $\R$ of `bounded variation': we can
always expand such an $f$ in the pointwise-convergent Fourier series $f(x)=
\sum_{n=0}^\infty a_n\, e(x)^n$.

There is a standard historical choice for this change-of-variables
function, namely
$$\eqalignno{j(\tau)\eqde&\, 1728\,{20\,G_4(\tau)^3\over 20\,G_4(\tau)^3-49\,G_6(\tau)^2}&\cr
=&\,
q^{-1}+744+196\,884\,q+21\,493\,760\,q^2+864\,299\,970\,q^3+\cdots
&(3)}$$
where as always $q=\exp[2\pi\i\,\tau]$. In fact, this choice is canonical,
apart from the arbitrary constant 744. This function $j$ is called the absolute
invariant or  {\it Hauptmodul} for SL$_2(\Z)$, or simply 
the {\it $j$-function}.

\smallskip In any case, one of the best studied functions
of classical number theory is the $j$-function. However, one of its most
remarkable properties was discovered only recently: McKay's approximations
 $196\,884\approx
196\,883$, $21\,493\,760\approx 21\,296\,876$, and $864\,299\,970\approx
842\,609\,326$. In fact,
$$\eqalignno{196\,884=&\,196\,883+1&(4a)\cr
21\,493\,760=&\, 21\,296\,876+196\,883+1&(4b)\cr
864\,299\,970=&\,842\,609\,326+ 21\,296\,876+2\cdot 196\,883+2\cdot 1&(4c)
\cr}$$
The numbers on the left sides of (4) are the first few coefficients of the
$j$-function (the number `744' in (3) is  of no
mathematical significance and can be ignored).
The numbers on the right are the dimensions
of the smallest irreducible representations of the {\it Monster finite simple
group} ${\Bbb M}$.

The {\it finite simple groups} are to group theory what the prime numbers
are to number theory --- in a sense they are the elementary building
blocks of all finite groups. They have been classified: the proof,
completed recently by a whole generation of group theorists, 
runs to approximately $15\,000$ pages and is spread over 500
individual papers. The resulting list
consists of 18 infinite families (e.g.\ the cyclic groups $\Z/p\Z$ of
prime order), together with 26 exceptional groups. The Monster $\M$ is 
the largest and richest of these exceptionals.

A {\it representation} of a group $G$ is the assignment of a matrix
$R(g)$ to each element $g$ of $G$ in such a way that the matrix
product respects the group product, i.e.\ $R(g)\,R(h)=R(gh)$. The
dimension of a representation is the size $n$ of its $n\times n$
matrices $R(g)$.

The equations (4) tell us that there is an
infinite-dimensional  graded representation 
$$V=V_{-1}\oplus V_1\oplus V_2\oplus V_3\oplus\cdots$$
of $\M$, where $V_{-1}=\rho_0$, $V_{1}=\rho_1\oplus\rho_0$, $V_2=\rho_2\oplus\rho_1
\oplus\rho_0$, $V_3=\rho_3\oplus\rho_2\oplus\rho_1\oplus\rho_1\oplus\rho_0
\oplus\rho_0$, etc, for the  irreducible representations $\rho_i$
of $\M$ (ordered by dimension), and that
$$j(\tau)-744={\rm dim}(V_{-1})\,q^{-1}+\sum_{i=1}^\infty {\rm dim}(V_i)
\,q^i$$
is the graded dimension of $V$.
John  Thompson then suggested that we `twist' this, i.e.\ that more generally we
consider the {\it McKay-Thompson series}
$$T_g(\tau)\eqde 
{\rm ch}_{V_{-1}}(g)\,q^{-1}+\sum_{i=1}^\infty {\rm ch}_{V_i
}(g)\,q^i\eqno(5)$$
for each element $g\in \M$. The character `ch$_R$' of a representation $R$
is given by `trace': ch$_R(g)={\rm Tr}(R(g))$. Up to equivalence (i.e.\
choice of basis), a representation $R$ can be recovered from the
character ch$_R$. The character however is much simpler --- e.g.\ the
smallest nontrivial representation of the Monster $\M$ is given by about
$10^{54}$ matrices, each of size $196\,883\times 196\,883$, while the
corresponding character is completely specified by 194 numbers (194 being the
number of conjugacy classes in $\M$).

The point of (5) is that, for any group representation
$\rho$, the character value ${\rm ch}_\rho
(id.)$ equals the dimension of $\rho$, and so
$T_{id.}(\tau)=j(\tau)-744$ and we recover (4) as special cases. But
there are many other possible choices of $g\in\M$, although conjugate
elements $g,hgh^{-1}$ will have identical character values and hence have
identical McKay-Thompson series $T_g=T_{hgh^{-1}}$. In fact there are precisely
171 {\it distinct} functions $T_g$. Perhaps these functions $T_g(\tau)$
might also be interesting.

Indeed, John Conway and Simon Norton [6] found that the first few terms of each McKay-Thompson series
$T_g$ coincided with the first few terms of certain special
functions, namely the `Hauptmoduls' of various
`genus-0 modular groups'.   Monstrous Moonshine --- which conjectured
that the McKay-Thompson series {\it were} those Hauptmoduls --- was officially born.

We should explain those terms.
We can generalise Definition 1 by replacing SL$_2(\Z)$ with any discrete
subgroup ${\cal G}$ of GL$_2(\Q)^+$, i.e.\ $2\times 2$ rational matrices with
positive determinant --- these act on $\H$ by fractional linear transformations
as usual. We can study the modular functions for
${\cal G}$ in the usual way, by studying the space ${\cal
G}\backslash \H$ of orbits. If ${\cal G}$ is not too big and not too
small,  then ${\cal G}\backslash{\H}$ will again be a compact Riemann surface
with finitely many points removed (corresponding as before to the cusps). When
this surface is a sphere, we call the modular group ${\cal G}$ {\it genus-0},
and the (appropriately normalised) change-of-variables function from
${\cal G}\backslash\H$ to the Riemann sphere $\C\cup\{\infty\}$ is
again called the {\it Hauptmodul} for ${\cal G}$. All modular
functions for a genus-0 group ${\cal G}$ will be rational functions of this
Hauptmodul. (On the other hand, when ${\cal G}$ is not genus-0, two generators are
needed, and unfortunately there is no canonical choice for them.)

The word `moonshine' here is English slang for `insubstantial or unreal'.
 It was chosen by Conway to convey as well the impression that 
things here are dimly lit, and that Conway-Norton were `distilling
 information illegally' from the Monster character table.

In hindsight, the first incarnation of Monstrous Moonshine goes back to
Andrew Ogg in 1975. He was in France discussing his result that the
primes $p$ for which the group ${\cal G}=\Gamma_0(p)+$ has genus 0, are 
$$p\in\{2,3,5,7,11,13, 17,19,23,29,31,41,47,59,71\}$$ 
$\Gamma_0(p)+$ is the group generated by all matrices
$\left(\matrix{a&b\cr c&d}\right)\in{\rm SL}_2(\Z)$ with $p$ dividing
the entry $c$, along with the matrix $\left(\matrix{0&1\cr -p&0}\right)$.
He also attended at that time a lecture
by Jacques Tits, who was describing a newly conjectured simple group.
When Tits wrote down the order 
$$\|\M\|=2^{46}\cdot 3^{20}\cdot 5^9\cdot 7^6\cdot 11^2\cdot 13^3\cdot 17
\cdot 19\cdot 23\cdot 29\cdot 31\cdot 41\cdot 47\cdot 59\cdot 71\approx
8\times 10^{53}$$
of that group, Ogg noticed its prime factors precisely equalled his
list of primes. Presumably as a joke, he offered a bottle of Jack Daniels'
whisky to the first person to explain the coincidence. Incidentally, we
now know that each of Ogg's groups $\Gamma_0(p)+$ is the genus-0
modular group for the function $T_g$, for some element $g\in\M$ of order $p$.

The significance of Monstrous Moonshine lies in its mysteriousness: it associates
various special modular functions to the Monster, even though mathematically
they seem fundamentally incommensurate. Now,
`understanding' something means to embed it naturally into a broader
context. Why is the sky blue? Because of the way light scatters in
gases. In order to understand Monstrous Moonshine, to resolve the mystery,
we should search for
similar phenomena, and fit them all into the same story.\smallskip

In actual fact, Moonshine (albeit non-Monstrous) really began long ago. Euler (and probably people
before) played with the power series $t(x)\eqde 1+2x+2x^4+2x^9+2x^{16}+\cdots$, primarily
because it can be used to express the number of ways a given number can be written
 as a sum of squares of integers. In his study of elliptic integrals,
 Jacobi noticed that if we change variables
by $x=e^{\pi \i\tau}$, then the resulting function $\theta_3(\tau)\eqde 1+
2e^{\pi\i\tau}+2e^{4\pi\i\tau}+\cdots$ behaves nicely with respect to
certain transformations of $\tau$ --- it's
a {modular form} for a certain subgroup of SL$_2(\Z)$. More generally, the same conclusion holds when
we sum not over the squares of $\Z$, but the norms of any
$n$-dimensional lattice $\L\subset \R^n$: 
the lattice theta series
$$\Theta_\L(\tau)\eqde\sum_{x\in\L}e^{\pi\i\, x\cdot x}$$
is also a modular form, provided  all norms $x\cdot x$ in $\L$ are rational.
See [7] for a fascinating and readable account of lattice lore.
(The role of the $n$-dimensional lattice $\L$ here is completely different
from that of the two-dimensional lattice $\L(\tau)$ in e.g.\ (1).)

And in the late 1960s Victor Kac [16] and Robert Moody [22]
independently (and for completely different reasons) defined
 a new class of infinite-dimensional Lie algebras. Within a decade
it was realised that the characters 
of the {\it affine} Kac-Moody algebras are (vector-valued) modular functions.

Indeed, McKay had also remarked in 1978 that similar coincidences to (4)
hold if $\M$ and
$j(\tau)$ respectively are replaced with the Lie group $E_8(\C)$ and
$$j(q)^{{1\over 3}}=q^{-{1\over 3}}\, (1+248q+4124q^2+34\,752q^3+\cdots)$$
In particular, $248={\rm dim}\,L(\L_7)$, $4124={\rm
dim}\,(L(\L_1)\oplus L(\L_7)\oplus L(0))$, $34\,752={\rm
dim}\,(L(\L_6)\oplus L(\L_1)\oplus 2\,L(\L_7)\oplus L(0))$, where $L(\la)$
denotes the representation of $E_8(\C)$ with `highest weight' $\la$. Incidentally,
$j^{{1\over 3}}$ is the Hauptmodul of the genus-0 group $\G(3)$, where
$$\G(N)\eqde \{A\in{\rm SL}_2(\Z)\,|\,A\equiv\left(\matrix{1&0\cr
0&1}\right)\ {\rm (mod}\ N)\}$$
In no time Kac [15] and Jim Lepowsky [20] independently remarked that the unique
level 1 highest-weight representation ${L}(\widehat{\L}_0)$ of the affine Kac-Moody algebra
$E_8^{(1)}$ has character $(qj(q))^{{1\over 3}}$. Since each graded
piece of any representation ${L}(\widehat{\la})$ of the affine
Kac-Moody algebra $X_\ell^{(1)}$ must
carry a representation of the associated finite-dimensional Lie group $X_\ell(\C)$, and the characters $\chi_{\hat{\la}}$
of an affine algebra are modular functions for some ${\cal
G}\subseteq{\rm SL}_2(\Z)$, this explained
McKay's $E_8$ observation. His Monster observations took longer to
clarify, because much of the mathematics was still to be developed. 

A {\it Lie algebra} ${\frak g}$ is a vector space with a bilinear vector-valued
product which is both anti-commutative and anti-associative. The
familiar cross-product in three-dimensions defines a Lie algebra,
called sl(3), and in fact this algebra can be used to generate all the
Kac-Moody algebras in a way encoded in the corresponding Coxeter-Dynkin diagram. A very readable introduction to Lie theory is the
book [5]. The standard references for Kac-Moody algebras are [17] and [18].
We'll return to Lie theory later in this paper.

We've known for years that lattices and affine Kac-Moody algebras are associated to
modular forms and functions. But
these observations, albeit now familiar, are also a little mysterious, we
should  confess. For instance, compare the unobvious fact that $\theta_3(-1/\tau)
=\sqrt{{\tau\over \i}}\,\theta_3(\tau)$, with the trivial observation (2) that
$G_k(-1/\tau)=\tau^kG_k(\tau)$ for the
 Eisenstein series $G_k$ in (1). The modularity of $\theta_3$, unlike that
of $G_k$, begs a {\it conceptual} explanation, even though its {\it logical}
explanation (i.e.\ proof) follows in a few moves from the Poisson summation
formula: 
$$\sum_{x\in \L}f(x)={1\over\sqrt{|\L|}}\sum_{y\in\L^*}\widehat{f}(y)\eqno(6)$$
where $\L\subset\R^n$ is any lattice (e.g.\ $\L=\Z\subset\R^1$), $\L^*\subset\R^n$ is
its dual lattice, $f$ is any `rapidly decreasing smooth function' on
$\R^n$, and $\widehat{f}(y)\eqde\int e^{-2\pi\i\,x\cdot y}f(x)\,dx$ is the
Fourier transform of $f$ (see e.g.\ \S6.1 of [26] for details). The key to the
simple  $\tau\mapsto -1/\tau$
transformation of $\theta_3$ is that the Fourier transform of the
Gaussian distribution  $e^{-\pi x^2}$ is itself. \medskip

{\narrower \noindent {\it At minimum, Moonshine  should be regarded as a
certain collection of related examples
where algebraic structures have been associated automorphic functions or forms.}

}

\bigskip\epsfysize=1.6in\centerline{ \epsffile{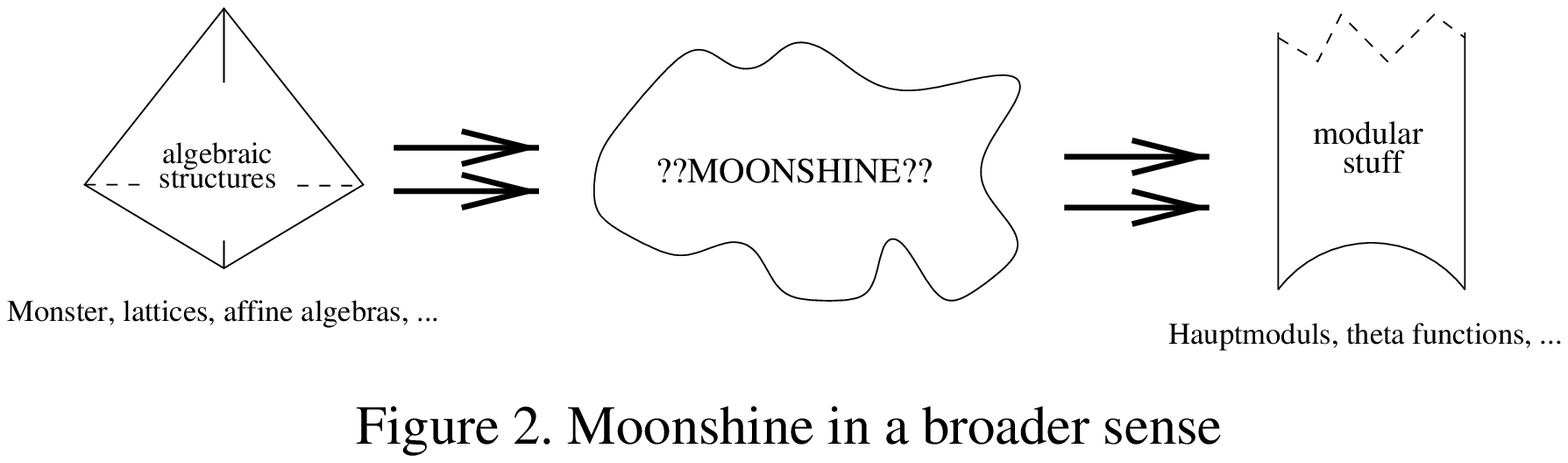}}\medskip

 From this larger perspective, illustrated in Figure 2,
 what is so special about the isolated example called
{\it Monstrous} Moonshine is that
the associated modular functions are of a special class (namely are
 Hauptmoduls). For lack of a better name, we call the theory of the
 blob of Figure 2, the {\it Theory of Generalised Moonshine}.

The first major step in the proof of Monstrous Moonshine was
accomplished in the mid 1980s with the construction by Frenkel-Lepowsky-Meurman
(see e.g.\ [10]) of the Moonshine module $V^{\natural}$ and its interpretation by
Borcherds [2] as a {\it vertex (operator) algebra}. A vertex operator algebra
is an infinite-dimensional vector space with infinitely many heavily constrained
 vector-valued bilinear products. Now, the Monster $\M$ is presumably  a
natural mathematical object, so we can expect that an elegant construction
for it would exist. Since $\M$ is the automorphism group of $V^\natural$,
and $V^\natural$ seems to be a natural though extremely intricate 
mathematical structure,  the hope it seems has been fullfilled.

In 1992 Borcherds [3] completed the proof of the Monstrous Moonshine
conjectures by showing that the graded characters $T_g$ of $V^\natural$
 are indeed the Hauptmoduls
conjectured by Conway and Norton, and hence that
$V^\natural$ is indeed the desired representation $V$ of $\M$
conjectured by McKay and Thompson.
The  explanation of Moonshine suggested by this picture
is given in Figure 3.
The algebraic structure can arise as the automorphism group of the associated
vertex operator algebra, or it can be hard-wired into the structure of the
vertex operator algebra. The modular forms/functions arise as the characters
of the (possibly twisted) modules of the vertex operator algebra.

\medskip\epsfysize=1.5in \centerline{\epsffile{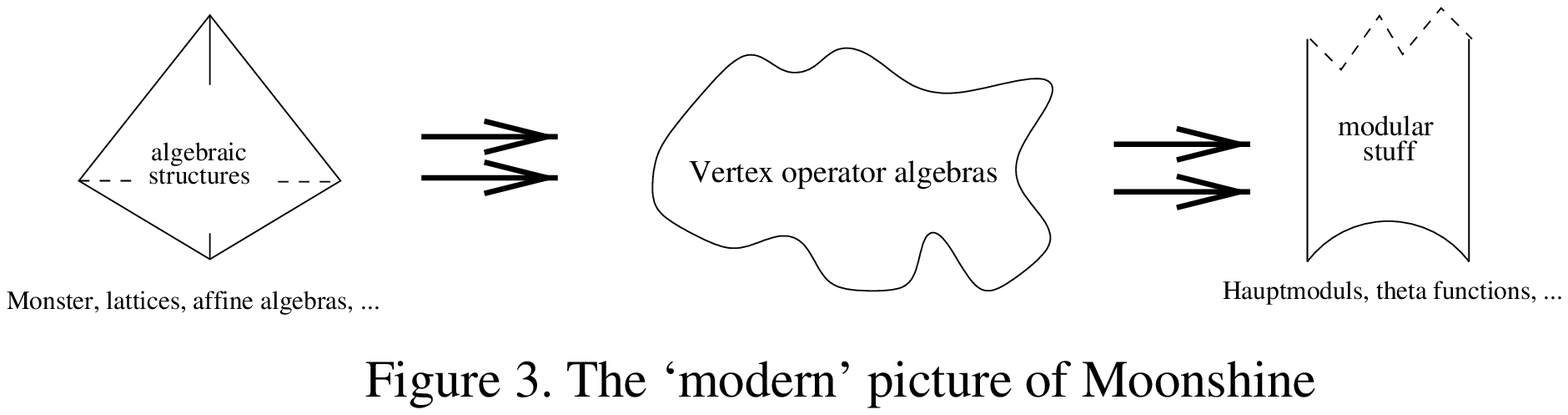}}\medskip

It must be emphasised that
 Figure 3 is meant to address Moonshine in the broader
sense of Figure 2, so certain special features of e.g.\ Monstrous Moonshine (in particular
that Hauptmoduls arise) will have to be treated by special arguments.

To see this genus-0
property of the $T_g$, Borcherds constructed a Kac-Moody-like Lie algebra
from $V^\natural$. The `(twisted) denominator identities' of this algebra
supply us with infinitely many equations which the coefficients $a_n(g)$
of the series $T_g$ must obey. For different reasons, the same equations
must be obeyed by the coefficients of the Hauptmoduls. These equations mean
that both the series $T_g$, and the Hauptmoduls, are uniquely determined
by their first few coefficients, so an easy computer check verifies
that each
$T_g$ equals the appropriate Hauptmodul. A more conceptual proof of this
Hauptmodul property was supplied in [8]: the denominator identities can
be reinterpreted as saying that the $T_g$ possess infinitely many `modular
equations'; it can be shown that any function obeying enough modular
equations must necessarily be a Hauptmodul.

Moonshine for other finite groups is explored in [24]. But what is so special
about the Monster $\M$, that its McKay-Thompson series $T_g$ are Hauptmoduls? It has
been conjectured [23] that it has to do with the `6-transposition property'
of $\M$. This thought has been further developed by Conway, Hsu,
Norton, and Parker
in their theory of quilts (see e.g.\ [14]). The genus-0 property for $\M$ has also been
related [27] to the conjectured uniqueness of the Moonshine module $V^\natural$.

Connections of Monstrous Moonshine with physics --- namely conformal field theory (CFT) [9] and
string theory ---
abound. A vertex operator algebra is an algebraic abstraction of
(one `chiral half' of) conformal field theory.
The Moonshine module $V^{\natural}$
can be interpreted as the string theory for a $\Z_2$-orbifold of free bosons compactified on
the torus $\R^{24}/\L_{24}$ associated to the Leech lattice $\L_{24}$. 
Many aspects of Monstrous Moonshine make complete sense within CFT, but some
(in particular the genus-0 property) remain more obscure.
In any case, although our story is primarily a mathematical one, most of the
chairs on which we sit were
warmed by physicists. In particular, what CFT (or what is essentially the
same thing, string theory) is, at least in part, is
a machine for producing modular functions. Figure 3 becomes Figure 4.
More precisely, the algebraic structure is an underlying
symmetry of the CFT, and its characters are the various modular functions.
The lattice theta functions come from bosonic strings living
on the torus $\R^n/\L$. The affine Kac-Moody characters arise in a
string theory where the string lives on a Lie group. And the Monster
is the automorphism group of a special `holomorphic' CFT intimately connected with $V^\natural$.

\medskip\epsfysize=1.5in \centerline{\epsffile{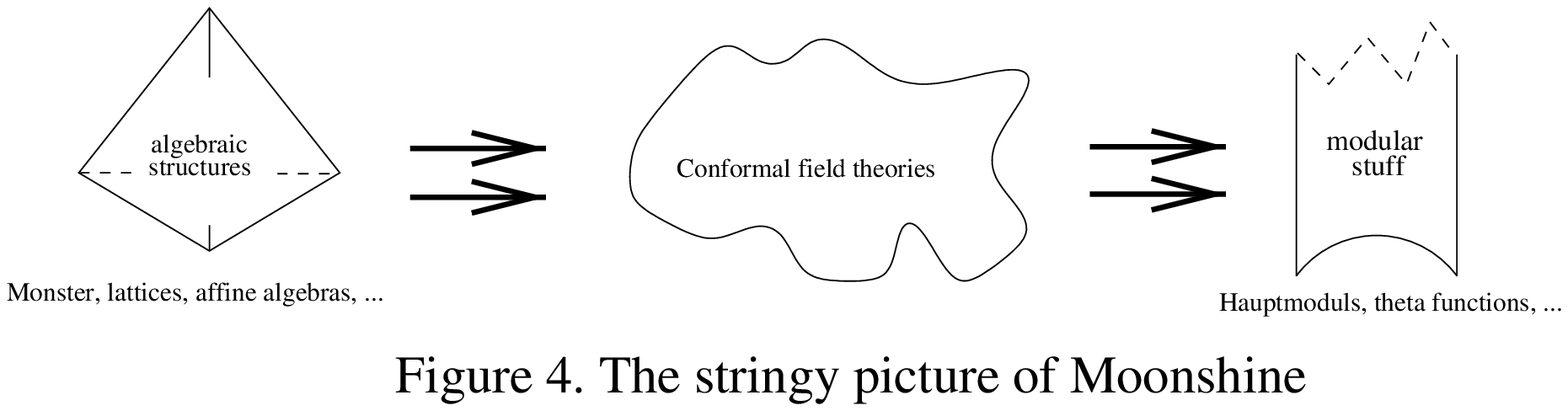}}\medskip

Historically speaking, Figure 4 preceded Figure 3.
The stringy picture is exciting because the CFT machine in
Figure 4 outputs much more
than merely modular functions --- it generates automorphic functions and forms for the
various mapping class groups of surfaces with punctures. And all this
is still poorly explored. We can thus expect more from Moonshine than
Figure 3 alone suggests.
On the other hand,
once again, Figure 4 by itself can only explain the broader
aspects of Moonshine.
More importantly, no one really knows what a CFT is (an influential but
incomplete attempt is by Graeme Segal [25]). Though that too may
be exciting to some physicists (and dismissed as inconsequential by others), most mathematicians find it a disturbing
flaw with Figure 4. Indeed, the definition by Borcherds and Frenkel-Lepowsky-Meurman
of a vertex operator algebra can be regarded as the first precise definition
of the {\it chiral algebra} of a CFT, and for this reason alone is a major achievement.

In spite of the work of Borcherds and others,
the special features of Monstrous Moonshine still beg questions.
The full conceptual relationship between the Monster and the
Hauptmoduls (like $j$) arguably remains `dimly lit', although much progress has been
realised. This is a subject where it is much easier to conjecture than
to prove, and we are still awash in unresolved conjectures. 

Nevertheless, Borcherds' paper [3] brings to a close the opening
chapter of the saga of Monstrous Moonshine. We are now in a period of
consolidation and synthesis, and it is in this spirit that this paper
is offered. 

So far, all of our `postcards' have been directly in the spirit of
Monstrous Moonshine. But the blob of Figure 2 is much more versatile
than that. We describe next three other postcards from the realm of
generalised Moonshine, which are in a sense {\it orthogonal} to
Monstrous Moonshine.\bigskip

Consider the following scenario. Let $A,B$ and $C$ be $n\times n$ Hermitian
matrices with eigenvalues $\alpha_1\ge\alpha_2\ge\cdots\ge \alpha_n$,
$\beta_1\ge\cdots\ge \beta_n$, $\gamma_1\ge\cdots\ge\gamma_n$. What are
the conditions on these eigenvalues such that $C=A+B$? The answer consists
of a number of inequalities involving the numbers $\alpha_i,\beta_j,\gamma_k$.
Discretise this problem, by requiring all $\alpha_i,\beta_j,\gamma_k$ to be nonnegative
integers. Then the following are equivalent (see e.g.\ [12]):

\medskip\item{(a)} Hermitian matrices $A,B$, and $C=A+B$ exist with
eigenvalues $\alpha,\beta,\gamma$, repsectively;

\smallskip\item{(b)} the GL$_n(\C)$ tensor product coefficient $T_{\alpha\beta}^\gamma$
is nonzero.\medskip

The finite-dimensional irreducible modules $L$ of
GL$_n(\C)$ are naturally labelled by such $n$-tuples $\alpha,\beta,\gamma$. The number
$T_{\alpha\beta}^\gamma$ is the number of times the module
$L(\gamma)$ appears in the tensor product $L(\alpha)\otimes L(\beta)$.

Now consider instead $n\times n$ unitary matrices with determinant 1. Any
such matrix $D\in {\rm SU}_n(\C)$ can be assigned a unique $n$-tuple $\delta=(\delta_1,\ldots,
\delta_n)$ as follows. Write its eigenvalues as $e^{2\pi\i\,\delta_i}$, where
$\delta_1\ge\cdots\ge\delta_n$, $\sum_{i=1}^n\delta_i=0$, and $\delta_1-
\delta_n\le 1$. Let $\Delta_n$ be the set of all such $n$-tuples $\delta$,
as $D$ runs through SU$_n(\C)$. Note that $D$ will have finite order iff all $\delta_i\in\Q$,
and that $D$ will be a scalar matrix $dI$ iff all differences $\delta_i-\delta_j\in
\Z$. Of course, a sum of Hermitian matrices corresponds here to a product
of unitary matrices.

Choose any {\it rational} $n$-tuples $\alpha,\beta,\gamma\in\Delta_n\cap
\Q^n$. Then the following are equivalent [1]:

\medskip\item{(i)} there exist matrices $A,B,C\in{\rm SU}_n(\C)$, where $C=AB$,
 with $n$-tuples $\alpha,\beta,\gamma$;

 \smallskip\item{(ii)} there is a positive integer $k$ such that
 all differences $k\alpha_i-k\alpha_j,k\beta_i-k\beta_j,k\gamma_i-k\gamma_j$
 are integers, and the sl$_{n}^{(1)}$ level $k$ fusion coefficient $N^{(k)\ k\gamma}_{
 k\alpha,k\beta}$ is nonzero.\medskip

`sl$_{n}^{(1)}$' is an affine Kac-Moody algebra. Here we interpret
$k\alpha$ etc as lying in the weight
lattice $A_{n-1}^*$, and so they correspond to the Dynkin labels $\la_i=k\alpha_i-k\alpha_{i+1}$,
etc of a level $k$ integrable highest-weight $\la$.  

The GL$_n(\C)$ tensor product coefficients $T_{\alpha\beta}^\gamma$ --- or {\it Littlewood-Richardson
coefficients} --- are classical quantities, appearing
in numerous and varied contexts. The sl$_n^{(1)}$ fusion coefficients
$N^{(k)\ \nu}_{\lambda\mu}$ are
equally fundamental, equally ubiquitous, but are more modern. For example,
they arise as tensor product coefficients for quantum groups at roots of
1, as dimensions of spaces of generalised theta functions, as dimensions
of conformal blocks in CFT, and as coefficients in the quantum cohomology
ring. They are perhaps the
most interesting example of a {\it fusion ring} (defined shortly). Fusion rings
are an aspect of generalised Moonshine complementary to Monstrous Moonshine, in the
sense that the fusion ring associated to Monstrous Moonshine is
trivial (i.e.\ one-dimensional).

\smallskip\noindent{{\smcap Definition 2.}}  A {\it fusion ring} [11,9,13]
$R$  is a commutative ring $R$ with identity 1, together with a finite
basis $\Phi$ (over $\Q$ say) containing 1, such that:

\smallskip \item{{\bf F1.}} The structure constants
$N_{ab}^c$ are all nonnegative integers;

\item{{\bf F2.}} There is a ring endomorphism $x\mapsto
x^*$ stabilising the basis $\Phi$;

\item{{\bf F3.}} $N_{ab}^1=\delta_{b,a^*}$.

\smallskip In addition, a self-duality condition identifying $R$ with
its dual should probably be imposed --- see [13] for details. As an
abstract ring it is not so interesting, as it is isomorphic (as an
algebra) to a direct sum of number fields. What is essential here is
the preferred basis $\Phi$.

The endomorphism $x\mapsto x^*$ can be shown to be an involution. We can
derive that there will be a unitary matrix $S$, with rows and columns
parametrised by $\Phi$, such that both $S_{1a},S_{a1}>0$ $\forall a$, and
$$N_{ab}^c=\sum_{i} {S_{ai}\,S_{bi}\,\overline{S_{ci}}\over S_{1i}}\eqno(7)$$
where  $\overline{S}$ denotes complex conjugate. 
The aforementioned self-duality condition amounts to a relation between $S$
and $S^t$ [13]. 

The fusion ring of a nontwisted affine algebra $X_\ell^{(1)}$ at
`level' $k\in\{1,2,3,\ldots\}$ is
$$R={\rm Ch}(X_\ell)/{\cal I}_k$$
where ${\rm Ch}(X_\ell)$ is the character ring of the Lie algebra $X_\ell$
(which has preferred basis given by the characters ch$_\lambda$, and whose
structure constants are the tensor product coefficients), and
where ${\cal I}_k$ is the ideal generated by all characters
of $X_\ell$ with level $k+1$.
(For $X_\ell=A_\ell$, the level of
representation $\lambda$ is given by $\sum_{i=1}^\ell\lambda_i$.)
The preferred basis for the fusion ring $R$ consists of
all characters ch$_\lambda$ with $\lambda$ of level $\le k$.
It is known that the $N_{\la\mu}^{\nu\ (k)}$ are nonnegative integers,
which increase with $k$ to the corresponding tensor product coefficient
$T_{\la\mu}^\nu$. Incidentally,
the {\it twisted} affine algebras also appear very naturally here, in
the context of `NIM-reps' or `fusion graphs', but this is another story.

What has a fusion ring to do with `modular stuff'? That is explained
in our next postcard: {\it modular data}.
\bigskip

Choose any even integer $n>0$. 
The matrix $S=({1\over \sqrt{n}}e^{2\pi\i\, mm'/n})_{0\le m,m'<n}$
 is the finite Fourier transform. Define the diagonal matrix $T$ by
$T_{mm}=e^{\pi\i m^2/n-\pi\i/12}$. The assignment 
$\left(\matrix{0&-1\cr 1&0}\right)\mapsto S$, $\left(\matrix{1&1\cr
0&1}\right)\mapsto T$ defines an $n$-dimensional representation $\rho$
of
SL$_2(\Z)$, since the matrices $\left(\matrix{0&-1\cr 1&0}\right)$ and
$\left(\matrix{1&1\cr 0&1}\right)$ generate SL$_2(\Z)$.
In fact this is essentially a Weil representation of SL$_2(\Z/n\Z)$. 
This is the simplest (and least interesting) example of what we'll 
call {\it modular data} --- a refinement of fusion rings to be defined shortly.
{Verlinde's formula} (7) here is the product rule for discrete exponentials, namely
$$e^{2\pi\i \,mm'/n}\cdot e^{2\pi\i\, mm''/n}=e^{2\pi\i\,
m\,(m'+m'')/n}$$

This representation is realised by modular functions. 
For each $m\in\{0,1,\ldots,n-1\}$, define
the functions
$$\psi_m(\tau)={1\over \eta(\tau)}\sum_{k=-\infty}^\infty q^{n\,(k+m/n)^2/2}$$
where as always $q=e^{2\pi\i\tau}$ and where $\eta(\tau)$ is the Dedekind
eta function:
$$\eta(\tau)\eqde \sum_{k=-\infty}^\infty (q^{6\,(k+{1\over 12})^2}-q^{6\,
(k+{5\over 12})^2})=\{{675\over 256\,\pi^{12}}(20\,G_4(\tau)^3-49\,G_6(\tau)^2)\}^{{1\over 24}}$$
 If we write $\L$ for the lattice
$\sqrt{n}\Z$, then $\L^*={1\over\sqrt{n}}\Z$ is the dual lattice, the number
$0\le m<n$ parametrises the cosets $\L^*/\L$, and $\psi_m$ is the theta series
of the $m$th coset. It's easy to see that $\psi_m(\tau+1)=T_{mm}\,\psi_m(\tau)$;
the Poisson summation formula (6) gives us
$\psi_m(-1/\tau)=\sum_{m'=0}^{n-1}
S_{mm'}\,\psi_{m'}(\tau)$.
Thus $\vec{\psi}=(\psi_0,\psi_1,\ldots,\psi_{n-1})^t$ is a `vector-valued
modular function with multiplier $\rho$' for SL$_2(\Z)$, in the sense that $\vec{\psi}(A.\tau)=\rho(A)\,
\vec{\psi}(\tau)$ for any $A\in{\rm SL}_2(\Z)$.

More generally, 
to various algebraic structures (in the above special case this is the lattice
$\L=\sqrt{n}\Z$) can be associated an SL$_2(\Z)$ representation. More interesting
examples come from affine Kac-Moody algebras and finite groups. The role
of $\psi_m$ is played by the characters of vertex operator algebras [29]
(or Kac-Moody algebras or CFT). Verlinde's formula (7)
associates a fusion ring to modular data. In Monstrous Moonshine, 
the modular data is trivial: each matrix $U$ is the $1\times 1$ matrix $(1)$.

\medskip\noindent{{\smcap Definition 3.}} Let $\Phi$ be a finite set of labels, one of which --- denoted `1'
and called the `identity' ---
is distinguished. By {\it modular data} we mean matrices
$S=(S_{ab})_{a,b\in \Phi}$, $T=(T_{ab})_{a,b\in \Phi}$ of complex
numbers  such that [13]:\smallskip

\item{{\bf M1.}} $S$ is unitary and symmetric, and $T$ is diagonal
  and of finite order: i.e.\ $T^N=I$ for some $N$;

\item{{\bf M2.}} $S_{1a}>0$ for all $a\in \Phi$;

\item{{\bf M3.}} $S^2=(ST)^3$;

\item{{\bf M4.}} The numbers $N_{ab}^c$ defined by Verlinde's formula (7)
are nonnegative integers.\medskip

Axiom {\bf M2} as stated is too strong, although Perron-Frobenius tells us that
some column of $S$ must be of constant phase. 
Modular data defines a representation of the modular group
SL$_2(\Z)$. Each entry $S_{ab}$ lies in some cyclotomic field extension
${\Bbb K}_n\eqde \Q[\exp(2\pi\i/n)]$. There is a simple and important action
of Gal$({\Bbb K}_n/\Q)\cong(\Z/n\Z)^*$ on $S$, which generalises the $g\mapsto
g^\ell$ symmetry of the character table of a finite group. In all known
examples, this SL$_2(\Z)$ representation is trivial on the principal
congruence subgroup $\Gamma(N)$ defined earlier, where $N$ is the order of $T$,
which means that the characters are modular
functions for $\Gamma(N)$, and that we really have  a representation for
the finite group SL$_2(\Z)/\Gamma(N)\cong{\rm SL}(\Z/N\Z)$.

\bigskip

A {\it knot} $K$ in $\R^n$ is a smooth one-to-one embedding of $S^1$ into $\R^n$.
The Jordan curve theorem states that all knots in $\R^2$ are trivial. Are
there any nontrivial knots in $\R^3$?

In Figures 5 and 6 we draw some knots in $\R^3$, by flattening them
into the plane of the paper. A moment's consideration will confirm that the
second knot of Figure 5 is indeed trivial. What about the trefoil?

\medskip\epsfysize=1.5in \centerline{\epsffile{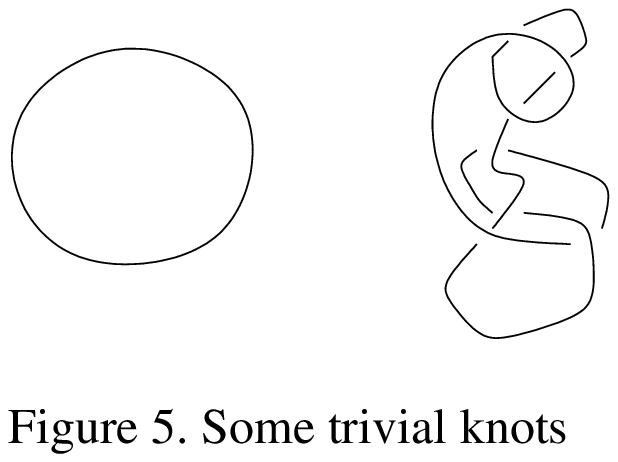}\qquad\qquad
\epsffile{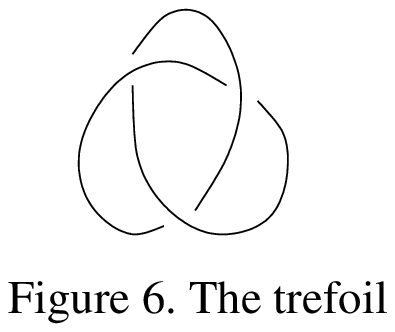}}\medskip

A knot diagram cuts the knotted $S^1$ into several connected components ({\it arcs}),
whose endpoints lie at the various {\it crossings} (double-points of
the projection).
By a {\it 3-colouring}, we mean to colour each arc 
in the knot diagram either red, blue or green, so that at each crossing
either 1 or 3 distinct colours are used.
 For example, the first two colourings in
Figure 7 are allowed, but the third one isn't.
By considering the `Reidemeister moves', which tell how to move between
equivalent knot diagrams, 
different diagrams for equivalent knots (such as the two in Figure 5) can be seen to have the same number of
distinct 3-colourings. Hence
 the number of different 3-colourings is a knot invariant.

\medskip\epsfysize=1.5in \centerline{\epsffile{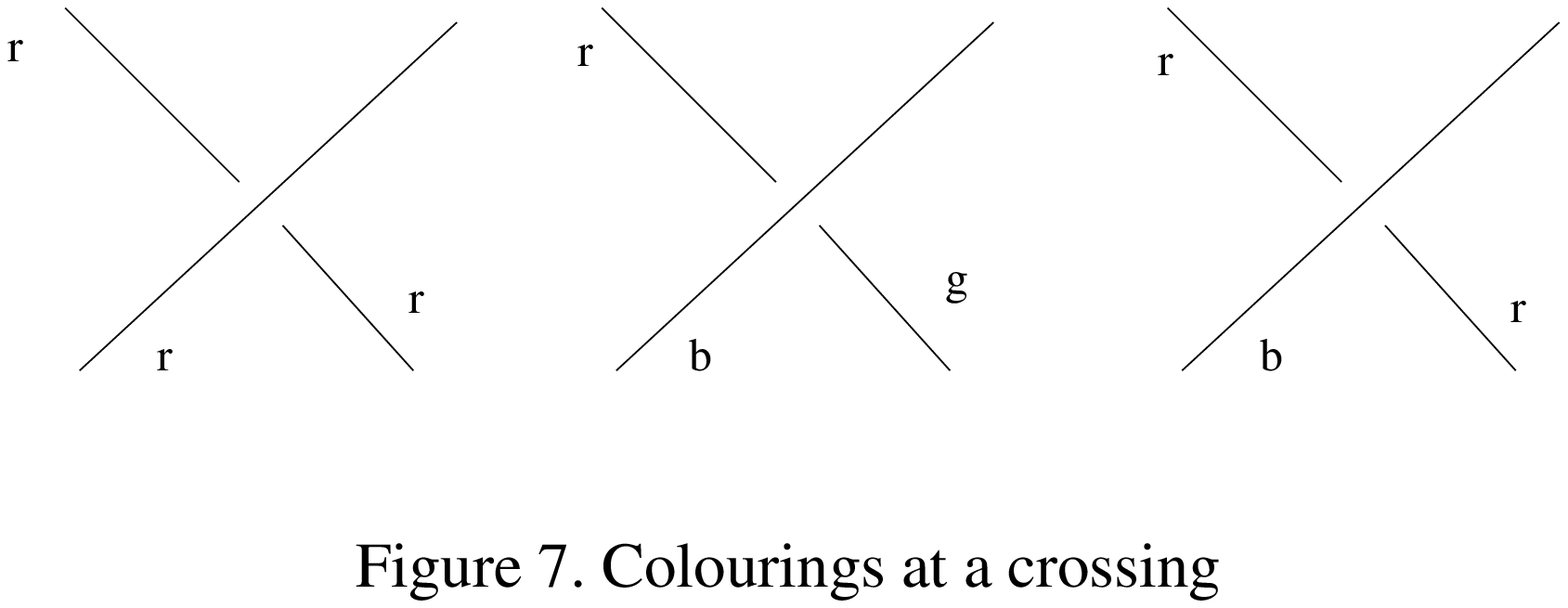}}\medskip

For example, consider the diagrams in Figure 5 for the trivial knot: the reader
can quickly verify that all arcs must be given the same colour, and thus there
are precisely three distinct 3-colourings. On the other hand, the trefoil  has nine
distinct 3-colourings --- the bottom two arcs of Figure 6 can be assigned arbitrary
colour, and that choice fixes the colour of the top arc. Thus  the
trefoil is nontrivial!

Essentially what we are doing here is counting the number of homomorphisms  $\varphi$
from the {\it knot group} $\pi_1(\R^3/K)$ to the symmetric group $S_3$.
The reason is that any knot diagram gives a presentation for $\pi_1(\R^3/K)$,
where there is a generator $g_i$ for each arc and a relation of the form
$g_i^{\pm 1}g_jg_i^{\mp 1}=g_k$ for each crossing. The map $\varphi$ is
defined using e.g.\ the identification $r\leftrightarrow (12)$, $b\leftrightarrow
(23)$, $g\leftrightarrow(13)$, and the above 3-colouring condition at each crossing
is equivalent to requiring that $\varphi$ obeys each group relation. Our
homomorphism $\varphi$ will be onto iff at least two different colours are
used.

By considering more general (nonabelian) colourings, the target ($S_3$ here)
can be made to be any  other group $G$, resulting in a different knot
invariant. This class of knot invariants is an example of one coming
 from {\it topological
field theory} (a refinement of modular data), in this case associated
to an arbitrary finite group $G$. Another deep and fascinating
source of topological field theories (and modular data etc) is
subfactor theory for von Neumann algebras --- a gentle introduction to
some aspects of this is [19]. 
The definition of topological field
theory is too long and
complicated to give here, but an excellent account is [28]. A standard introduction to
knot theory is [4].

What has topological field theory to do with modular stuff? The matrix
$S$ comes from the knot invariants attached to the so-called Hopf link
(two linked circles in $\R^3$). The knots and links here are really
`framed', i.e.\ are ribbons, and the diagonal matrix $T$ describes
what happens when the ribbon is twisted. If $S$ and $T$ constitute
modular data (defined earlier), then the topological field theory will yield knot
invariants in any 3-manifold (via the process called surgery). The fusion coefficients come from three
parallel circles  $p_i\times S^1$ in the 3-manifold $S^2\times
S^1$. There is no canonical choice of characters (modular functions) though which
realise this SL$_2(\Z)$ representation.

For instance, returning to the topological field theory and modular
data associated to
finite group $S_3$, we have $T={\rm
diag}(1,1,1,1,e^{2\pi\i/3},e^{-2\pi\i/3},1,-1)$, and
$$S={1\over 6}\left(\matrix{1&1&2&2&2&2&3&3\cr 1&1&2&2&2&2&-3&-3\cr 2&2&4&-2&-2&-2&0&0\cr 
2&2&-2&4&-2&-2&0&0\cr 2&2&-2&-2&-2&4&0&0\cr 2&2&-2&-2&4&-2&0&0\cr 3&-3&0&0&0&0&3&-3\cr
3&-3&0&0&0&0&-3&3}\right)$$
\smallskip

Like moonlight itself, Moonshine is an indirect phenomenon. Just as in the theory
of moonlight one must introduce the sun, so in the theory of Moonshine
one should go beyond the Monster. Much as a review paper discussing moonlight may
include a few paragraphs on sunsets or comet tails, so have we sent postcards
of fusion
rings, SL$_2(\Z)$ representations, and knot invariants.

\bigskip\noindent{{\smcap Bibliography}} \medskip

\item{[1]} S.\ Agnihotri and C.\ Woodward, ``Eigenvalues of products
of unitary matrices and quantum Schubert calculus'', {\it Math.\ Res.\
Lett.} {\bf 5} (1998), 817--836.

\item{[2]} R.\ E.\ Borcherds, ``Vertex algebras, Kac-Moody algebras, and the
Monster'', {\it Proc.\ Natl.\ Acad.\ Sci.\ (USA)} {\bf 83} (1986), 3068--3071.

\item{[3]} {R.\ E.\ Borcherds},
``Monstrous moonshine and monstrous Lie superalgebras'',
{\it Invent.\ math.}\ {\bf 109} (1992) 405--444.

\item{[4]} G.\ Burke and H.\ Zieschang, {\it Knots} (de Gruyter, Berlin, 1995).

\item{[5]} R.\ Carter, G.\ Segal and I.\ M.\ Macdonald, {\it Lectures
on Lie Groups and Lie Algebras} (Cambridge University Press, Cambridge, 1995).

\item{[6]} {J.\ H.\ Conway and S.\ P.\ Norton},
``Monstrous moonshine'',
{\it Bull.\ London Math.\ Soc.}\ {\bf 11} (1979), 308--339.

\item{[7]} {J.\ H.\ Conway and N.\ J.\ A.\ Sloane}, {\it Sphere
Packings, Lattices and Groups}, 3rd edn (Springer, Berlin, 1999).

\item{[8]} {C.\ J.\ Cummins and T.\ Gannon}, ``Modular equations and
the genus zero property'',
{\it Invent.\ math.} {\bf 129} (1997), 413--443.

\item{[9]} P.\ Di Francesco, P.\ Mathieu and D.\ S\'en\'echal, {\it
Conformal Field Theory} (Springer, New York, 1996).

\item{[10]} {I.\ Frenkel, J.\ Lepowsky, and A.\ Meurman}, {\it
Vertex Operator Algebras and the Monster} (Academic Press, San Diego,
1988).

\item{[11]} {J.\ Fuchs}, ``Fusion rules in conformal field theory'',
{\it Fortsch.\ Phys.} {\bf 42} (1994), 1--48.

\item{[12]} W.\ Fulton, ``Eigenvalues, invariant factors, highest weights,
and Schubert calculus'', math.AG/9908012.

\item{[13]} T.\ Gannon, ``Modular data: the algebraic combinatorics of
conformal field theory'', math.QA/0103044.

\item{[14]} T.\ Hsu, {\it Quilts: Central Extensions, Braid Actions,
and Finite Groups}, Lecture Notes in Math 1731 (Springer, Berlin, 2000).

\item{[15]} V.\ G.\ Kac, ``Simple irreducible graded Lie algebras of
finite growth'', {\it Math.\ USSR--Izv} {\bf 2} (1968), 1271--1311. 

\item{[16]} V.\ G.\ Kac,  
``An elucidation of: Infinite-dimensional algebras, Dedekind's $\eta
$-function,  classical M\"obius function and the very strange
formula. $E_{8}^{(1)}$ and the cube root of the modular invariant $j$'', 
{\it Adv. in Math} {\bf 35} (1980), 264--273.

\item{[17]} {V.\ G.\ Kac}, {\it Infinite Dimensional Lie Algebras}, 
3rd edn (Cambridge University Press, Cambridge, 1990).

\item{[18]} S.\ Kass, R.\ V.\ Moody, J.\ Patera and R.\ Slansky, {\it
Affine Lie Algebras, Weight Multiplicities, and Branching Rules}, Vol.\ 1
(University of California Press, Berkeley, 1990).

\item{[19]} V.\ Kodiyalam and V.\ S.\ Sunder, {\it Topological Quantum
Field Theories from Subfactors} (Chapman \& Hall, New York, 2001).

\item{[20]} J.\ Lepowsky, 
``Euclidean Lie algebras and the modular function $j$'', In:
{\it Proc. Sympos. Pure Math.} {\bf  37} 
(Amer. Math. Soc., Providence, 1980),  pp. 567--570.

\item{[21]} H.\ McKean and V.\ Moll, {\it Elliptic Curves: Function Theory,
Geometry, Arithmetic} (Cambridge University Press, Cambridge, 1999).

\item{[22]} {R.\ V.\ Moody}, ``A new class of Lie algebras'',
{\it J.\ Algebra} {\bf 10} (1968), 211--230.

\item{[23]} S.\ P.\ Norton, ``Generalized moonshine'', In: {\it
Proc.\ Symp.\ Pure Math} {\bf 47} (Amer. Math. Soc., Providence,
1987), pp.208--209.

\item{[24]} L.\ Queen, ``Modular functions arising from some
finite groups'', {\it Math.\ of Comput.}\ {\bf 37} (1981), 547--580.

\item{[25]} G.\ Segal, ``Geometric aspects of quantum field theories'', In:
{\it Proc.\ Intern.\ Congr.\ Math., Kyoto} (Springer, Hong Kong, 1991), pp.1387--1396.

\item{[26]} J.-P.\ Serre, {\it A Course in Arithmetic} (Springer, Berlin, 1973).

\item{[27]} M.\ Tuite, ``On the relationship between
Monstrous moonshine and the uniqueness of the Moonshine module'', {\it Commun.\
Math.\ Phys.}\  {\bf 166} (1995), 495--532.

\item{[28]} V.\ G.\ Turaev, {\it Quantum Invariants of Knots and 3-Manifolds}
(de Gruyter, Berlin, 1994).

\item{[29]} {Y.\ Zhu}, ``Modular invariance of characters
of vertex operator algebras'', {\it J.\ Amer.\ Math.\ Soc.} {\bf 9}
(1996), 237--302.

\end